\documentclass[a4paper]{article}
\newcommand{\bI}{{\bf I}}
\newcommand{\bJ}{{\bf J}}

\newcommand{\bV}{{\bf V}}
\newcommand{\bW}{{\bf W}}
\newcommand{\bT}{{\bf T}}
\newcommand{\bE}{{\bf E}}

\newcommand{\ba}{{\bf a}}
\newcommand{\bb}{{\bf b}}
\newcommand{\bc}{{\bf c}}
\newcommand{\bx}{{\bf x}}

\newcommand{\ad}{{\rm{}^{ad}}}

\newcommand{\Pro}{{\cal P}}

\newcommand{\RR}{\mbox{\rm I\hspace{-0.166 em}R}}

\newcommand{\beweisende}{~\rule{1,2ex}{1,2ex}}

\newcommand{\vertsp}{\vspace{1 em}}

\newtheorem{theo}{Theorem}

\begin{document}
\vertsp\vertsp\vertsp\vertsp

\centerline{ON THE MATRICES OF CENTRAL LINEAR MAPPINGS}

\vertsp\vertsp

\centerline{{\sc Hans Havlicek}, Wien}
\vertsp
\vertsp\vertsp

   \begin{small}
   {\em Summary.} We show that a central linear mapping of a projectively
   embedded Euclidean $n$-space onto a projectively embedded Euclidean
   $m$-space is decomposable into a central projection followed by a
   similarity if, and only if, the least singular value of a certain matrix
   has multiplicity $\ge 2m-n+1$. This matrix is arising, by a simple
   manipulation, from a matrix describing the given mapping in terms of
   homogeneous Cartesian coordinates.

   \vertsp


   \vertsp\vertsp


   \end{small}

\section{Introduction}\label{sect-INTRO}

A linear mapping between projectively embedded Euclidean spaces is called
{\em central}, if its exceptional subspace is not at infinity. Such a linear
mapping is in general not decomposable into a central projection followed by a
similarity. Necessary and sufficient conditions for the existence of such a
decomposition have been given in \cite{Have1} for arbitrary finite dimensions;
cf. also \cite{Brau1}, \cite{Brau2}, \cite{Brau3}. However, those results do
not seem to be immediately applicable on a {\em central axonometry}, i.e., a
central linear mapping given via an axonometric figure. On the other hand, in
a series of recent papers \cite{Pauk1}, \cite{Szab1}, \cite{Szab2} this
problem of decomposition has been discussed for central axonometries of the
Euclidean 3-space onto the Euclidean plane from an elementary point of view%
   \footnote{A lot of further references can be found in the quoted papers.}%
.

Loosely speaking, the concept of central axonometry is a geometric
equivalent to the algebraic concept of a {\em coordinate matrix} for a linear
mapping of the underlying vector spaces. However, from the results in
\cite{Brau2} and \cite{Have1} it is also not immediate whether or not a given
matrix  describes (in terms of homogeneous Cartesian coordinates) a mapping
that permits the above-mentioned factorization. The aim of this communication
is to give a criterion for this.

\vertsp

Let $\bI$, $\bJ$ be finite-dimensional Euclidean vector spaces. Given a
linear mapping $f : \bI \to \bJ$ denote by $f\ad:\bJ\to\bI$ its adjoint
mapping. Then $f\ad \circ f$ is
self-adjoint with eigenvalues
   \begin{displaymath}
   v_1 \ge \cdots \ge v_r > v_{r+1} = \cdots = v_n = 0.
   \end{displaymath}
Here $r$ equals the rank of $f$ and $n=\dim \bI$. Moreover, each eigenvalue
is written down repeatedly according to its multiplicity%
   \footnote{For a self-adjoint mapping the algebraic and geometric
   multiplicities of an eigenvalue are identical. Hence we may unambiguously
   use the term `multiplicity'.}%
. The positive real numbers $\sqrt{v_1},\ldots,\sqrt{v_r}$ are frequently
called the {\em singular values} of $f$. The multiplicity of a singular value
of $f$ is defined via the multiplicity of the corresponding eigenvalue of
$f\ad \circ f$. It is immediate from the singular value decomposition that $f$
and $f\ad$ share the same singular values (counted with their multiplicities).
See, e.g., \cite{Stra}.

These results hold true, mutatis mutandis, when replacing $f$ by any real
matrix, say $A$, and $f\ad$ by the transpose matrix $A^{\rm T}$.

\section{Decompositions}\label{sect-KOO_FREI}
When discussing central linear mappings it will be convenient to consider
Euclidean spaces embedded in projective spaces. Thus let $\bV$ be an
$(n+1)$-dimensional real vector space ($3\le n < \infty$) and $\bI$ one of its
hyperplanes. Assume, furthermore, that $\bI$ is equipped with a positive
definite inner product ($\cdot$) so that $\bI$ is a Euclidean vector space. In
the projective space on $\bV$, denoted by $\Pro(\bV)$, we consider the
projective hyperplane $\Pro(\bI)$ as the hyperplane at infinity. The absolute
polarity in $\Pro(\bI)$ is determined by the inner product on $\bI$. Hence
$\Pro(\bV) \setminus \Pro(\bI)$ is a projectively embedded Euclidean space%
   \footnote{We do not endow this space with a unit segment.}%
. Similarly, let $\Pro(\bW) \setminus \Pro(\bJ)$ be an $m$-dimensional
projectively embedded Euclidean space $(2\le m < n < \infty)$.
Given a linear mapping
   \begin{equation}\label{f}
   f : \bV \to \bW
   \end{equation}
of vector spaces then the associate (projective) linear mapping
   \begin{equation}\label{phi}
   \phi : \Pro(\bV) \setminus\Pro(\ker f)\to \Pro(\bW) \mbox{, } \RR\bx
   \mapsto \RR(f(\bx))
   \end{equation}
has $\Pro(\ker f)$ as its exceptional subspace. In the sequel we shall assume
that
   \begin{equation}\label{cen-sur}
\ker f \not\subset \bI \quad\mbox{and}\quad f(\bV) = \bW,
   \end{equation}
or, in other words, that $\phi$ is central and surjective%
   \footnote{This  assumption  of  surjectivity  is  made  `without  loss   of
   generality' in most papers on this subject. It will, however, be  essential
   several times in this paper.}%
. Obviously, (\ref{cen-sur}) is equivalent to
   \begin{equation}\label{combi}
   f(\bI) = \bW.
   \end{equation}

We recall some results \cite{Brau2}, \cite{Have1}: If $\bT$ is any
complementary subspace of $\ker f$ in $\bV$, then denote by
   \begin{equation}
   \psi_{\bT} : \Pro(\bV) \setminus\Pro(\ker f) \to \Pro(\bT)
   \end{equation}
the projection with the exceptional subspace $\Pro(\ker f)$ onto $\Pro(\bT)$.
The restricted mapping
   \begin{equation}
   \phi_{\bT} := \phi \mid \Pro(\bT) : \Pro(\bT) \to \Pro(\bW)
   \end{equation}
is a collineation and
   \begin{equation}
   \phi = \phi_{\bT}\circ\psi_{\bT};
   \end{equation}
every decomposition of $\phi$ into a projection and a collineation is of this
form. In the Euclidean vector space $\bI$ we have the distinguished subspace
   \begin{equation}
   \bE := f^{-1}(\bJ)\cap \bI.
  \end{equation}
Write
   \begin{equation}
   f_\bE : \bE \to \bJ \mbox{, } \bx\mapsto f(\bx);
   \end{equation}
this $f_\bE$ is well-defined and surjective, since $\bE\subset f^{-1}(\bJ)$
and $\ker f\not\subset \bE$.
The subspace $\bT$ can be chosen with $\phi_\bT$ being a similarity if, and
only if, the least singular value of $f_\bE$ has multiplicity%
   \footnote{In \cite[Satz 10]{Brau2} this multiplicity is printed
   incorrectly as $2m-n-1$.}
$\ge 2m-n+1$.

Next, we assume that $\Pro(\bT)\not\subset \Pro(\bI)$ is orthogonal to
$\Pro(\ker f)$. This means that $(\bT\cap\bI)^\bot \subset \ker f\cap\bI$ or
$(\bT\cap\bI)^\bot \supset \ker f\cap\bI$. Hence $\psi_\bT$ is an {\em
orthogonal central projection}%
   \footnote{The central projections used in elementary descriptive geometry
   are trivial examples of orthogonal central projections.}%
. It is easily seen from \cite{Brau2} that $\phi$ permits a decomposition
into an orthogonal central projection followed by a similarity if, and only
if, all singular values of $f_\bE$ are equal.

Finally, we are going to show that the crucial properties of $f_\bE$ can be
read off from another mapping: Denote by
   \begin{equation}
   p : \bI \to \bE
   \end{equation}
the orthogonal projection with the kernel $\bE^\bot\subset\bI$. Then
   \begin{equation}\label{f-p}
   (f_\bE \circ p) \circ (f_\bE \circ p)\ad =
   f_\bE \circ p \circ p\ad \circ (f_\bE)\ad =
   f_\bE \circ (f_\bE)\ad,
   \end{equation}
since $p\ad$ is the natural embedding $\bE \to \bI$. Thus, by (\ref{f-p}) and
the results stated in Section \ref{sect-INTRO}, $f_\bE$ and $(f_\bE \circ
p)\ad$ have the same singular values (counted with their multiplicities).
Hence, by the surjectivity of $f_\bE$ and (\ref{f-p}), all singular values of
$f_\bE$ are equal if, and only if, there exists a real number $v>0$ such that
   \begin{equation}\label{ortho}
   (f_\bE \circ p) \circ (f_\bE \circ p)\ad = v\,\mbox{id}_\bJ.
   \end{equation}
We shall use this in the next section.

\section{A matrix characterization}\label{sect-KOORD}

Introducing homogeneous Cartesian coordinates in $\Pro(\bV)$ is equivalent to
choosing a basis $\{\bb_0,\ldots,\bb_n\}$ of $\bV$ such that
$\{\bb_1,\ldots,\bb_n\}\subset\bI$ is an orthonormal system. The origin is
given by $\RR\bb_0$ and the unit points are $\RR(\bb_0+\bb_1), \ldots,
\RR(\bb_0+\bb_n)$. In the same manner we are introducing homogeneous Cartesian
coordinates in $\Pro(\bW)$ via a basis $\{\bc_0,\ldots,\bc_m\}$.

   \begin{theo}
   Suppose that $f:\bV\to\bW$ is inducing a surjective central linear mapping
   $\phi$ according to formula (\ref{phi}). Let
      \begin{equation}
      A = \left(
         \begin{array}{lcl}
         a_{00}      &  \cdots   &  a_{0n}   \\
         {\rm\vdots} &  {}       &  {\rm\vdots}   \\
         a_{m0}      &  \cdots   &  a_{mn}
         \end{array}
      \right)
      \end{equation}
   be the coordinate matrix of $f$ with respect to bases of $\bV$ and $\bW$
   that are yielding homogeneous Cartesian coordinates. Write
      \begin{equation}\label{a_i}
      \ba_i := (a_{i1},\ldots,a_{in})\in \RR^n
      \mbox{ for all } i=0, \ldots, m
      \end{equation}
   and
      \begin{equation}
      \widetilde{A} := \left(
         \begin{array}{c}
         \ba_1- \frac{\ba_0\cdot\ba_1}{\ba_0\cdot\ba_0}\ba_0 \\
         {\rm\vdots}\\
         \ba_m - \frac{\ba_0\cdot\ba_m}{\ba_0\cdot\ba_0}\ba_0
         \end{array}
      \right).
      \end{equation}
   Then the following assertions hold true:
      \begin{enumerate}
      \item $\phi$ is decomposable into a central projection followed by a
      similarity if, and only if, the least singular value of the matrix
      $\widetilde{A}$ has multiplicity $\ge 2m-n+1$.
      \item $\phi$ is decomposable into an orthogonal central projection
      followed by a similarity if, and only if, there exists a real number
      $v>0$ such that
         \begin{equation}
         \widetilde{A} \widetilde{A}^{\rm T} = \mbox{\rm diag}\,(v,\ldots,v).
         \end{equation}
      \end{enumerate}
   \end{theo}
{\em Proof}. We read off from the top row of $A$ that
   \begin{displaymath}
   a_{00}x_0 + \cdots + a_{0n}x_n = 0
   \end{displaymath}
is an equation of $f^{-1}(\bJ)\neq \bI$ so that $\ba_0\cdot\ba_0\neq 0$. Write
$\widetilde{f}:\bI\to\bJ$ for the linear mapping whose coordinate matrix with
respect to $\{\bb_1,\ldots,\bb_n\}$ and $\{\bc_1,\ldots,\bc_m\}$ equals
$\widetilde{A}$. A straightforward calculation shows that
   \begin{displaymath}
   \widetilde{f}(\bx) = f(\bx) \mbox{ for all } \bx\in \bE
   \end{displaymath}
and
   \begin{displaymath}
   \widetilde{f}(a_{01}\bb_1+\cdots +a_{0n}\bb_n) = 0,
   \end{displaymath}
i.e., $\bE^{\bot}\subset\ker\widetilde{f}$. Thus $\widetilde{f}$ equals
the mapping $f_\bE \circ p$ discussed above. Now the proof is completed  by
translating formulae (\ref{f-p}) and (\ref{ortho}) into the language of
matrices.\beweisende

\vertsp

We remark that (\ref{cen-sur}) and the linear independence of
$\ba_1,\ldots,\ba_m$ are equivalent conditions.

In contrast to the results in \cite{Pauk1}, \cite{Szab1}, \cite{Szab2}, the
$\phi$-image of the origin $\RR\bb_0$ does not appear in our characterization.
On the other hand, we have
   \begin{displaymath}
   f(\bE^\bot) = \RR((\ba_0\cdot\ba_0)\bc_0 +\cdots +(\ba_0\cdot\ba_m)\bc_m).
   \end{displaymath}
In projective terms this 1-dimensional subspace of $\bW$ gives the {\em
principal point} of the mapping $\phi$. Exactly if the principal point of
$\phi$ equals the origin $\RR\bc_0$, then $\widetilde{A}$ arises from $A$
merely by deleting the top row and the leading column.

   \small
{\em Author's address}: Hans Havlicek, Abteilung f\"ur Lineare Algebra und
Geometrie, Technische Universit\"at, Wiedner Hauptstra{\ss}e 8--10, A--1040
Wien, Austria.

   \noindent
EMAIL: havlicek@geometrie.tuwien.ac.at
\end{document}